\documentclass[11pt,reqno]{amsart}
\usepackage{amsmath,amsopn,amssymb,amsthm,multicol}
\usepackage{hyperref}
\usepackage[all]{xy}
\usepackage{graphicx}
\usepackage{graphics}

\DeclareMathOperator{\Ad}{Ad}

\DeclareMathOperator{\Aut}{Aut}

\DeclareMathOperator{\GL}{GL}

\DeclareMathOperator{\Ric}{Ric}

\newcommand{\fr}{\mathfrak}
\newcommand{\al}{\alpha}

\newcommand{\bb}{\mathbb}

\DeclareMathOperator{\SO}{SO}

\DeclareMathOperator{\Sp}{Sp}
\DeclareMathOperator{\SU}{SU}
\DeclareMathOperator{\U}{U}
\DeclareMathOperator{\G}{G}
\DeclareMathOperator{\F}{F}

 \newtheorem{lemma} {Lemma}[section]
\newtheorem{theorem}[lemma]{Theorem} 
 
\newtheorem{prop} [lemma]{Proposition}  
\newtheorem{definition}[lemma] {Definition} 
\newtheorem{corol}[lemma] {Corollary} 
\newtheorem{example}[lemma] {Example}

\title{Invariant metrics on homogeneous spaces with equivalent isotropy summands} 
\author{Marina Statha}
\address{University of Patras, Department of Mathematics, GR-26500 Rion, Greece}
\email{statha@master.math.upatras.gr}
\medskip

\begin{document}

\begin{abstract}
The space of $G$-invariant metrics on a homogeneous space $G/H$ is in one-to-one correspondence with  the set of inner products on the tangent space $\fr{m}\cong T_{{\it o}}(G/H)$, which are invariant under the isotropy representation.  When all the isotropy summands are inequivalent to each other, then the metric is called diagonal.  We will describe a special class of $G$-invariant metrics 
in the case where the isotropy representation of $G/H$ contains some equivalent isotropy summands.  Even though this problem has been considered sporadically in the bibliography, in the present article we provide a more systematic and organized description of such metrics.  This will enable us to simplify the problem of finding $G$-invariant Einstein metrics for  homogeneous spaces.  We also provide some applications.
\medskip


\end{abstract}

\maketitle

\section{Introduction}

A homogeneous manifold $M$ is a manifold which admits a transitive group of diffeomorphisms.  However, in general there might be  several distinct transitive groups, i.e. non conjugate transitive subgroups of the diffeomorphism group of $M$, and these subgroups can be abstractly isomorphic.  If we fix a compact Lie group $G$ acting on a homogeneous manifold $M$, then after choosing a basepoint, we can write $M$ as the coset space $G/H$, where $H$ is the isotropy group at the basepoint.  From the theorem of Myers and Steenrod \cite{MySt} it follows that the isometry group ${\rm Iso}(M)$ of $M$, is a Lie group and that the isotropy subgroup $H$ is a closed compact subgroup of ${\rm Iso}(M)$.  One of the fundamental properties of a homogeneous space is that, if we know the value of a geometrical quantity at a given point, then we can calculate its value at any other point of $G/H$ by using translation maps.  Hence all calculations reduce to a single point which, for simplicity, can be chosen to be the identity coset ${\it o} = eH\in G/H$.  

A Riemannian manifold $(M, g)$ is called Einstein if the metric $g$ satisfies the condition $\Ric(g) = \lambda g$ for some $\lambda\in \bb{R}$.  We refer to \cite{Be} and \cite{Wa1}, \cite{Wa2} for old and new results on homogeneoous Einstein manifolds.  The structure of the set of invariant Einstein metrics on a given homogeneous space is still not very well understood in general.  The situation is only clear for few classes of homogeneous spaces.  
For an arbitrary compact homogeneous space $G/H$ it is not clear if the set of invariant Einstein metrics (up to isometry and
up to scaling) is finite or not.  A finiteness conjecture states that this set is in fact finite if the isotropy representation of $G/H$ consists of pairwise inequivalent irreducible components (\cite{BWZ}).  

A large class of homogeneous spaces are the \textit{reductive} homogeneous spaces.  For these spaces there exists a subspace $\fr{m}$ of $\fr{g}$ such that $\fr{g} = \fr{h}\oplus\fr{m}$ and $\Ad(H)\fr{m}\subset\fr{m}$.  The tangent space of $M$ at ${\it o}$ is canonically identified with $\fr{m}$.  A major class of reductive homogeneous spaces are the \textit{isotropy irreducible homogeneous spaces}.  These spaces have been studied by J. Wolf in \cite{Wo}, where he proved that \textit{if $G/H$ is an isotropy irreducible homogeneous space, then $G/H$ admits a unique (up to scalar) $G$-invariant metric}, which is also Einstein.  Later, M. Wang and W. Ziller in \cite{WZ1} and \cite{WZ2}, gave a complete classification of such spaces.  The most important examples of isotropy irreducible homogeneous spaces are the irreducible symmetric spaces, classified by E. Cartan in 1926.  More generally, in \cite{Be} it is shown that a non compact irreducible homogeneous space is symmetric.  If the reductive homogeneous space is not isotropy irreducible, then its isotropy representation splits into a direct sum of irreducible subrepresentations.  Examples of such spaces are the generalized flag manifolds, Wallach spaces, the projective space $\bb{C}{\rm P}^{2n+1}$ and the Stiefel manifolds.  

Generalized flag manifolds with two and four isotropy summands are classified using a method based on Riemannian submersions by A. Arvanitoyeorgos and I. Chrysikos in \cite{ArCh1}, \cite{ArCh2}.  In general, homogeneous spaces with two irreducible isotropy summands were classified by W. Dickinson and M. Kerr in \cite{DiKe}.  This classification is achieved under the assumptions that $G$ is a compact, connected and simple Lie group, $H$ is a closed subgroup of $G$ and $G/H$ is simply connected.  It should be noted that in this classification there is only one example of a homogeneous space having equivalent subrepresentations, namely the space $\SO(8)/\G_{2} \cong S^{7}\times S^{7}$.  The $G$-invariant Einstein metrics on this space as well as on the homogeneous spaces ${\rm Spin}(7)/\U(3)\cong S^{7}\times S^{6}$, ${\rm Spin}(8)/\U(3)\cong S^{7}\times \G_{2}^{+}(\bb{R}^{8})$ and on the Stiefel manifold $V_{2}\bb{R}^{n+1}\cong \SO(n+1)/\SO(n-1)$, where the isotropy representation splits into equivalent subrepresentations, were classified by M. Kerr in \cite{Ke}.  The Allof-Wallach spaces $W_{k, l} = \SU(3)/\SO(2)$ when $(k, l) = (1, 0)$ and $(k, l) = (1, 1)$ are two examples of homogeneous spaces with equivalent subrepresentations.  In general, the space of invariant Riemannian metrics on $W_{k, l}$, is parametrized by four positive parameters.  For $(k, l) = (1, 0)$ and $(k, l) = (1, 1)$ this space depends on $6$ and $10$ positive real numbers, respectively.  By using the variational approach Yu. Nikonorov in \cite{Ni1} proved that there are at most two invariant Einstein metrics on $W_{1, 1}$.  Morever, he constructed a new invariant Einstein metric on $W_{1,0}$ which is not diagonal with respect to the $\Ad(T)$-invariant decomposition of $\SU(3)$, where $T$ is a maximal torus in $\SU(3)$.  

Finally, A. Arvanitoyeorgos, Yu. Nikonorov and V. V. Dzhepko proved that for $s>1$ and $k>l\geq 3$ the Stiefel manifold $\SO(sk + l)/\SO(l)$ admits at least four $\SO(sk + l)$-invariant Einstein metrics, two of which are Jensen's metrics.  The special case $\SO(2k + l)/\SO(l)$ admitting at least four $\SO(2k + l)$-invariant Einstein metrics was treated in \cite{ArDzN1}.  Corresponding results for the quaternionic Stiefel manifolds $\Sp(sk + l)/\Sp(l)$ were obtained in \cite{ArDzN2}.  Recently, it was proved by A. Arvanitoyeorgos, Y. Sakane and the author in \cite{ArSaSt}, that the Stiefel manifold $V_{4}\bb{R}^{n}\cong\SO(n)/\SO(n-4)$ admits two more $\SO(n)$-invariant Einstein metrics and that $V_{5}\bb{R}^{7}\cong\SO(7)/\SO(5)$ admits four more $\SO(7)$-invariant Einstein metrics.

In the present paper we study $G$-invariant metrics on homogeneous spaces $G/H$ for which the isotropy representation contains equivalent subrepresentations or isotropy summands.  For such spaces the diagonal metrics are not unique.  We odserve that the normalizer $N_{G}(H)$ acts on the space of all $G$-invariant metrics $\mathcal{M}^{G}$ by isometries, and we can choose a subgroup $K$ of $N_{G}(H)$ such that the action of $K$ on $\mathcal{M}^{G}$ determines a subset of $\mathcal{M}^{G}$.  Our approach is analysized in Section 3 and is summarized in the following Theorem:
\begin{theorem}
Let $M= G/H$ be a homogeneous space of a compact semisimple Lie group $G$ and let $K$ be a closed subgroup of $G$ such that $H\subset K\subset N_G(H)$, where $N_G(H)$ is the normalizer of $H$ in $G$.
\begin{itemize}
\item[$(1)$]  The non trivial action $(\varphi, A)\mapsto \varphi\circ A\circ \varphi ^{-1}$  of the set $\Phi_K = \{\varphi = \Ad(k)|_{\fr{m}} : k\in K\} \subset \Phi =\{\phi = \Ad(n)|_{\fr{m}} : n\in N_{G}(H)\subset \Aut(\fr{m})$, on the set $\mathcal{M}^G$ of all $G$-invariant metrics on $G/H$ is well defined.
\item[$(2)$]  The set $(\mathcal{M}^G)^{\Phi _K} = \{A\in \mathcal{M}^{G} : \varphi\circ A\varphi^{-1} = A \ \mbox{for all} \ \varphi\in \Phi_{K}\}$ of fixed points of the action in $(1)$  determines a subset of all $\Ad(H)$-invariant inner products on $\fr{m}$, called $\Ad(K)$-invariant inner products.  This set in turn, determines a subset $\mathcal{M}^{G, K}$ of $\mathcal{M}^G$.
\end{itemize}
\end{theorem} 
\noindent
Theorems of the above type are useful for the study of geometrical problems (e.g. finding $G$-invariant Einstein metrics) on homogeneous space whose isotropy representation contains equivalent summands (see for example \cite{ArSaSt}).

The paper is organized as follows: In Section 2 we recall some useful results from representation theory.  In Section 3 we analyse the action of the normalizer $N_{G}(H)$ on the set of $\mathcal{M}^{G}$ of all $G$-invariant metrics on $G/H$.  By restricting this action to a closed subgroup $K$ of $G$ such that $H\subset K\subset N_{G}(H)$, we obtain a subset $\mathcal{M}^{G, K}$ of all $G$-invariant metrics $\mathcal{M}^{G}$.  As a consequence, various geometrical objects (such as Ricci tensor) are easier to by described.  In Section 4 we relate such a choice of subgroup $K$ (i.e. $H\subset K \subset N_{G}(H)$) to Riemannian submersions $K/H \to G/H \to G/K$. 

\medskip

\noindent {\bf Acknowledgements.}
The work was supported by Grant $\# E.037$ from the Research Committee of the University of Patras (Programme K. Karatheodori).  The author gratefully acknowledges  very useful discussions and influence by professors Yusuke Sakane and Yurii Nikonorov.  She also expresses her gratitude to professor Andreas Arvanitoyeorgos for his constant guidence  and support.

\section{Review of representation theory}

A finite dimensional (real or complex) representation of a Lie group $G$ is a homomorphism $\varphi : G \to \Aut(V)$, where $V$ is a finite dimensional (real or complex) vector space.  The dimension of the representation is the dimension of the vector space $V$.  If there is no non trivial subspace $W\subset V$ with $\varphi(W)\subset W$ then the representation $\varphi$ called irreducible.  The complexification of a real representation $\varphi : G \to \Aut(V)$ is defined as the complex representation $\varphi\otimes\bb{C} : G \to \Aut(V\otimes\bb{C})$.

\begin{definition}
Two representations $\varphi_{1} : G \to \Aut(V_{1})$ and $\varphi_{2} : G \to \Aut(V_{2})$ are called equivalent $(\varphi_{1} \cong \varphi_{2}\ \mbox{or}\ V_{1}\cong V_{2})$ if $V_{1}$ and $V_{2}$ are $G$-isomorphic, i.e. there exists a linear isomorphism $f : V_{1} \to V_{2}$ such that $f(\varphi_{1}(g)v) = \varphi_{2}(g)f(v)$, for all $g\in G$ and $v\in V_{1}$.  Such an $f$ is also called $G$-equivariant map (or intertwining map).
\end{definition}

A useful observation is the following.

\begin{theorem}{\bf(Schur's Lemma)}
If $\varphi : G \to \Aut(V)$ is an irreducible complex representation and $f\in {\rm Hom}(V, V)$ is a $G$-equivariant map, then $f = c{\rm Id}$ for some $c\in\bb{C}$.
\end{theorem} 

For every representation $\varphi : G \to \Aut(V)$ of a compact topological group $G$ there exists a $G$-invariant inner product $\langle\cdot, \cdot\rangle$ on $V$, i.e. $\langle \varphi(g)u, \varphi(g)v \rangle =$ $\langle u, v \rangle$, for all $g\in G$ and $u, v\in V$.  From this it follows that any representation of a compact topological group is a direct sum of irreducible representations i.e.
$
\varphi \cong \varphi_{1}\oplus\cdots\oplus\varphi_{n} : G\to \Aut(V_{1}\oplus\cdots\oplus V_{n}),
$
where each of $\varphi_{i} : G \to \Aut(V_{i})$ $(i = 1,2,\ldots,n)$ is irreducible.  

If $\varphi$ is a real (resp. complex) irreducible representation and $\langle\cdot,\cdot\rangle_{1}$, $\langle\cdot,\cdot\rangle_{2}$ are two $G$-invariant inner products (resp. hermitian inner products) on $V$, then from the above theorem it follows that $\langle\cdot,\cdot\rangle_{1} = c \langle\cdot,\cdot\rangle_{2},$ for some $c\in\bb{R}$ (resp. $c\in\bb{C}$).  Therefore, if $\varphi \cong \varphi_{1}\oplus\cdots\oplus\varphi_{n}$ and assuming that $\varphi_{i}$ are mutually inequivalent, then all $G$-invariant inner products on $V$ are given by
$$
\langle\cdot,\cdot\rangle = x_{1}\, \langle\cdot,\cdot\rangle |_{V_{1}} + \cdots + x_{n}\, \langle\cdot,\cdot\rangle |_{V_{n}}, \ \ x_{i}\in\bb{R}, \ i = 1,\ldots,n
$$
where $\langle V_{i}, V_{j}\rangle = 0$ for $i\neq j$.  Any other $G$-invariant inner product on $V$ can be expressed as $(\cdot, \cdot) = \langle A\cdot, \cdot\rangle$, where $A : V \to V$ is a positive definite, symmetric, $G$-equivariant linear map.  

If $\varphi_{i}$ and $\varphi_{j}$ are equivalent for some $i$ and $j$, then the above inner product is not unique, and $\langle V_{i}, V_{j}\rangle$ does not necessarily vanish, thus the matrix of the operator $A$ has some non zero non diagonal elements.  To find the number of non diagonal elements, we need to determine the dimension of the space of intertwining maps \textit{between the pairs of equivalent representations}.  For example, let  $\varphi_{1}\cong\varphi_{2}$ and $\varphi_{i},\ i=1,2$ be irreducible as real representations.  The complexification of $\varphi_{1}$ is not necessarily irreducible.  After complexifying $\varphi_{1}$, there are three possibilities \textnormal{(\cite{Ke})}:
\begin{itemize}
\item[1.] If $\varphi_{1}\otimes\bb{C}$ is irreducible, we call $\varphi_{1}$ orthogonal. 
\item[2.] If $\varphi_{1}\otimes\bb{C} = \psi\oplus\bar{\psi}$ and $\psi$ is not equivalent to $\bar{\psi}$, we call $\varphi_{1}$ unitary. 
\item[3.] If $\varphi_{1}\otimes\bb{C} = \psi\oplus\bar{\psi}$ and $\psi$ is equivalent to $\bar{\psi}$, we call $\varphi_{1}$ symplectic.  
\end{itemize}
The space of intertwining maps is $1$-dimensional in the orthogonal case, $2$-dimensional in the
unitary case, and $4$-dimensional in the symplectic case.  Thus if in the decomposition of $V = V_{1}\oplus V_{2}\oplus\cdots\oplus V_{n}$ we have $r$ equivalent summands (or modules), then the number of non diagonal elements in the orthogonal case is $\frac{r(r-1)}{2}$, in the unitary case it is $r(r-1)$ and in the symplectic case it is $2r(r-1)$.  In the present article we describe a special class of $G$-invariant metrics on a homogeneous spaces $G/H$ which contain equivalent isotropy summands.  

\begin{definition}
The adjoint representation of $G$ is the homomorphism $\Ad \equiv \Ad^{G} : G\to \Aut(\fr{g})$ given by $\Ad(g) = (dI_{g})_{e}$, where $I_{g} : G\to G, x\mapsto g x g^{-1}$, and $\fr{g}$ is the Lie algebra of $G$.
\end{definition}

Denote by $\tilde{\lambda}_{n}$ the standard representation of ${\rm GL}_{n}\bb{R}$ and by $\lambda_{n}$ the standard representation of $\SO(n)$ (or ${\rm O}(n)$).  It is $\lambda_{n} = \tilde{\lambda}_{n}|_{\SO(n)} : {\GL}_{n}\bb{R}\to \Aut(\bb{R}^{n})$.  Then the adjoint representation $\Ad^{\SO(n)}$ of $\SO(n)$ (or ${\rm O}(n)$) is equivalent to $\wedge^{2}\lambda_{n}$, where $\wedge^{2}$ denotes the second exterior power of $\lambda_{n}$.  Also, we have that $\Ad^{\U(n)}\otimes\bb{C} = \mu_{n}\otimes\bar{\mu}_{n}$ and $\Ad^{\Sp(n)}\otimes\bb{C} = S^{2}\nu_{n}$, where $\mu_{n} = \tilde{\mu}_{n}|_{\U(n)} : \GL_{n}\bb{C}\to \Aut(\bb{C}^{n})$, \ $\nu_{n} = \tilde{\nu}_{n}|_{\Sp(n)} : \GL_{n}\bb{H}\to \Aut(\bb{H}^{n})$.  Here $\tilde{\mu}_{n}$, $\tilde{\nu}_{n}$ are the standard representations of ${\rm GL}_{n}\bb{C}$ and ${\rm GL}_{n}\bb{H}$ respectively, and $S^{2}$ is the second symmetric power of $\nu_{n}$.  Recall that if $\pi : G \to \Aut(V),$ \ $\pi' : G' \to \Aut(W)$ are two representations of $G$ and $G'$ respectively, then the following identities are valid: 
$$
\wedge^{2}(\pi\oplus\pi') = \wedge^{2}\pi\oplus\wedge^{2}\pi'\oplus(\pi\otimes\pi'), \ 
S^{2}(\pi\oplus\pi') = S^{2}\pi\oplus S^{2}\pi'\oplus(\pi\otimes\pi').
$$

Let $M$ be a smooth manifold and let $G$ be a Lie group acting on $M$ on the left by the map $\al : G\times M \to M, (g, m)\mapsto \al(g, m) = g m$.  For all $g\in G$, let $\al_{g} : M \to M$ be the corresponding diffeomorphism of $M$.  If $H = \{g\in G : gp = p\}$ is the isotropy subgroup at the point $p\in M$, then the \textit{isotropy representation} of $H$ at $p$ is the homomorphism 
\begin{eqnarray}\label{isotr}
\theta : H &\longrightarrow& \Aut(T_{p}M) \nonumber\\
h &\longmapsto& (d\al_{h})_{p} : T_{p}M \to T_{p}M,
\end{eqnarray} 
where $T_{p}M$ is the tangent space of $M$ at the point $p$.  In the case where the above action is also transitive, i.e. for $p, q\in M$ there exists $g\in G$ such that $q = gp$, then $M$ is diffeomorphic to the homogeneous space $G/H$, where $H$ is the isotropy subgroup at the identity coset ${\it o} = e H$.  By (\ref{isotr}) the \textit{isotropy representation} of $G/H$ is the homomorphism
\begin{eqnarray*}
\Ad^{G/H} : H &\longrightarrow& \Aut(T_{{\it o}}(G/H)) \nonumber\\
h &\longmapsto& (d\tau_{h})_{{\it o}} : T_{{\it o}}(G/H) \to T_{{\it o}}(G/H),
\end{eqnarray*}
where $\tau_{h} : G/H \to G/H,$\ $gH\mapsto hgH$.  A large class of homogeneous spaces are the \textit{reductive homogeneous spaces}.  For such spaces there exists a subspace $\fr{m}$ of the Lie algebra $\fr{g}$ such that $\fr{g} = \fr{h}\oplus\fr{m}$ and $\Ad(h)\fr{m}\subset \fr{m}$ for all $h\in H$, that is $\fr{m}$ is $\Ad(H)$-invariant.  If the subgroup $H$ is compact such decomposition always exists.  Then we have a canonical isomorphism $\fr{m}\cong T_{{\it o}}(G/H)$ given by $X \leftrightarrow X_{{\it o}}^{*} = \displaystyle{\frac{d}{dt}(\exp (tX)){\it o}|_{t = 0}}$, where $\exp (tX)$ is the one parameter subgroup of $G$ generated by $X$.

The next proposition is useful to compute the isotropy representation of the reductive homogeneous space (\cite{Arv}). 

\begin{prop}\label{isotrepr}
Let $G/H$ be a reductive homogeneous space and let $\fr{g} = \fr{h}\oplus\fr{m}$ be a reductive decomposition of $\fr{g}$.  Let $h\in H$, $X\in \fr{h}$ and $Y\in\fr{m}$.  Then
$$
\Ad^{G}(h)(X + Y) = \Ad^{G}(h)X + \Ad^{G}(h)Y
$$
that is, the restriction $\Ad^{G}\big|_{H}$ splits into the sum $\Ad^{H}\oplus\Ad^{G/H}$. 
\end{prop}

We give some examples of computations.
\begin{example}\label{pisotr}
\textnormal{We consider the homogeneous space $G/H = \SO(k_1 + k_2 + k_3)/(\SO(k_1)\times\SO(k_2)\times\SO(k_3))$ with $k_{1}, k_{2}, k_{3}\geq 2$, which is an example of a \textit{generalized Wallach space} (\cite{NiRoSl}).  These spaces were recently classified independently by Yu. Nikonorov in \cite{Ni2} and Z. Chen, Y. Kang, K. Liang in \cite{ChKaLi}.
Let $\sigma_{i} : \SO(k_1)\times\SO(k_2)\times\SO(k_3) \to \SO(k_i)$ be the projection onto the factor $\SO(k_{i}),$ $(i = 1,2,3)$ and let $p_{k_{i}} = \lambda_{k_{i}}\circ\sigma_{i}$.  Then we have the following:
\begin{eqnarray*}
\left.\Ad ^G\right |_H &=& \wedge ^2\lambda _{k_1+k_2+k_3}\big| _H = \wedge^2(p_{k_{1}}\oplus p_{k_{2}} \oplus p_{k_{3}}) = \wedge^{2}p_{k_{1}}\oplus\wedge^{2}p_{k_{2}} \nonumber\\ 
&& \oplus\wedge^{2}p_{k_{3}}\oplus(p_{k_{1}}\otimes p_{k_{2}})\oplus(p_{k_{1}}\otimes p_{k_{3}})\oplus(p_{k_{2}}\otimes p_{k_{3}}).
\end{eqnarray*} 
Observe that the dimension of the representation $\wedge^{2}p_{k_{1}}\oplus\wedge^{2}p_{k_{2}}\oplus\wedge^{2}p_{k_{3}}$ is $\binom{k_1}{2} + \binom{k_2}{2} + \binom{k_3}{2}$, which is equal to the dimension of the adjoint representation of  $H = \SO(k_1)\times\SO(k_2)\times\SO(k_3)$, $\Ad^{H} : \SO(k_1)\times\SO(k_2)\times\SO(k_3) \to \Aut(\fr{so}(k_{1})\oplus\fr{so}(k_{2})\oplus\fr{so}(k_{3}))$.  Therefore, the isotropy representation of $G/H$ is given by
\begin{equation}\label{isotrwall}
\Ad^{G/H}\cong (p_{k_{1}}\otimes p_{k_{2}})\oplus(p_{k_{1}}\otimes p_{k_{3}})\oplus(p_{k_{2}}\otimes p_{k_{3}}), 
\end{equation}
which is a direct sum of irreducible and non equivalent subrepresentations of dimensions $k_{i}k_{j}$, $i\neq j$.  The tangent space $\fr{m} $ of $G/H$ decomposes into three $\Ad(H)$-invariant submodules 
$
\fr{m} = \fr{m}_{12}\oplus  \fr{m}_{13}\oplus  \fr{m}_{23}. 
$}

\textnormal{Let us consider the case where $H_{1} = \SO(l_{1})\times\SO(l_{2})\times\SO(l_{3})$ and $l_{1} + l_{2} + l_{3} < k_{1} + k_{2} + k_{3} -1$.  Then we see that the isotropy representation of the homogeneous space $G/H_{1}$ contains some equivalent subrepresentations.  Indeed,
\begin{eqnarray*}
\left.\Ad ^G\right |_{H_{1}} &=& \wedge ^2\lambda _{l_1+l_2+l_3}\big| _{H_{1}} = \wedge^2(p_{l_{1}}\oplus p_{l_{2}}\oplus p_{l_{3}} \oplus 1_{n})
= \wedge^{2}p_{l_{1}}\oplus\wedge^{2}p_{l_{2}}\\
&& \oplus\wedge^{2}p_{l_{3}}\oplus\wedge^{2}1_{n}\oplus(p_{l_{1}}\otimes p_{l_{2}})\oplus(p_{l_{1}}\otimes p_{l_{3}})\oplus(p_{l_{2}}\otimes p_{l_{3}})\\
&&\oplus(p_{l_{1}}\otimes 1_{n})\oplus(p_{l_{2}}\otimes 1_{n})\oplus(p_{l_{3}}\otimes 1_{n})
\end{eqnarray*}
\begin{eqnarray*}
&=& \wedge^{2}p_{l_{1}}\oplus\wedge^{2}p_{l_{2}}\oplus\wedge^{2}p_{l_{3}}\oplus\underbrace{1\oplus\cdots\oplus 1}_{\binom{n}{2}}\oplus(p_{l_{1}}\otimes p_{l_{2}})\oplus(p_{l_{1}}\otimes p_{l_{3}})\\
&&\oplus(p_{l_{2}}\otimes p_{l_{3}}) \oplus\underbrace{p_{l_{1}}\oplus\cdots\oplus p_{l_{1}}}_{n}  \oplus\underbrace{p_{l_{2}}\oplus\cdots\oplus p_{l_{2}}}_{n}\oplus\underbrace{p_{l_{3}}\oplus\cdots\oplus p_{l_{3}}}_{n}.
\end{eqnarray*}
where $n = (k_1 + k_2 + k_3) - (l_1 + l_2 + l_3)$ and $1_{n} = \underbrace{1\oplus\cdots\oplus 1}_{n-\mbox{times}}$.  As before, the representation $\wedge^{2}p_{l_{1}}\oplus\wedge^{2}p_{l_{2}}\oplus\wedge^{2}p_{l_{3}}$ is the adjoint representation of $H_{1} = \SO(l_{1})\times\SO(l_{2})\times\SO(l_{3})$, thus the isotropy representation of the homogeneous space $G/H_{1}$ is 
\begin{eqnarray*}
\Ad ^{G/H_{1}} &=& 1\oplus\cdots\oplus 1\oplus(p_{l_{1}}\otimes p_{l_{2}})\oplus(p_{l_{1}}\otimes p_{l_{3}})\oplus(p_{l_{2}}\otimes p_{l_{3}})\\
&&\oplus p_{l_{1}}\oplus\cdots\oplus p_{l_{1}}\oplus p_{l_{2}}\oplus\cdots\oplus p_{l_{2}}\oplus p_{l_{3}}\oplus\cdots\oplus p_{l_{3}}.
\end{eqnarray*}
Observe that the last $3n$ representations of dimensions $l_{i}, (i = 1, 2, 3)$ are equivalent.  Thus the tangent space of $G/H_{1}$ decomposes into a sum of $\binom{n}{2} + 3n + 3$ $\Ad(H_{1})$-invariant submodules $\fr{m}_{i}$.  Similar result is true if we take $H_{2} = \SO(m_{1})\times\SO(m_{2})$ with $m_{1} + m_{2} < k_{1} + k_{2} + k_{3} - 1$, or $H_{3} = \SO(d)$ with $d < k_{1} + k_{2} + k_{3} - 1$.  In the special case where $H_{4} = \SO(k_{3})$, then the homogeneous space $G/H_{4}$ is the Stiefel manifold $V_{k_{1} + k_{2}}\bb{R}^{k_1 + k_2 + k_3}$.  In this case the isotropy representation is given as follows:
\begin{eqnarray*}
\left.\Ad ^G\right |_{H_{4}} &=& \wedge ^2\lambda _{k_1+k_2+k_3}\big| _{H_{4}} = \wedge^{2}(\lambda_{k_{3}}\oplus 1_{k_{1}+k_{2}})\\
&=& \wedge^{2}\lambda_{k_{3}} \oplus \wedge^{2}1_{k_{1} + k_{2}} \oplus (\lambda_{k_{3}}\oplus 1_{k_{1}+k_{2}})\\
&=& \wedge^{2}\lambda_{k_{3}} \oplus\underbrace{1\oplus\cdots\oplus 1}_{\binom{k_{1} + k_{2}}{2}}\oplus\underbrace{\lambda_{k_{3}}\oplus\cdots\oplus\lambda_{k_{3}}}_{k_{1} + k_{2}}\\
&=& \Ad^{\SO(k_{3})} \oplus 1\oplus\cdots\oplus 1 \oplus\lambda_{k_{3}}\oplus\cdots\oplus\lambda_{k_{3}},
\end{eqnarray*}   
hence the isotropy representation is
$
\Ad^{G/H_{4}} = 1\oplus\cdots\oplus 1 \oplus\lambda_{k_{3}}\oplus\cdots\oplus\lambda_{k_{3}},
$
where the last $k_{1} + k_{2}$ representations are equivalent.  
Analogous results can be obtained for $G = \SU(k_1 + k_2 + k_3)$ or $\Sp(k_{1} + k_{2} + k_3)$.  We summarize the above computations in the following table:
\begin{center}
\begin{tabular}{|l|lll|l|}
\hline\hline
$H$\ $\mbox{subgroup of}$\ $G$   &  $\fr{m} = \bigoplus_{i=1}^{s}\fr{m}_{i}$  &    $\mbox{non equiv.rep.}$   &$\mbox{equiv.rep.}$\\
\hline\hline
$\SO(k_{1})\times\SO(k_2)\times\SO(k_3)$    &           &                      &  \\
$k_{1}, k_{2}, k_{3}\geq 2$  &   $s = 3$  &      \ \ \ \    $\checkmark$ &\\
\hline
$\SO(l_{1})\times\SO(l_2)\times\SO(l_3)$    &                   &    &\\
$l_1 + l_2 + l_3 <$      &                  &       &\\
$k_1 + k_2 + k_3 -1$ &                  &       &\\
$n = (k_1 + k_2 + k_3)- $  &   $s = {\binom{n}{2}}+$               &       &\\
$ (l_1 + l_2 + l_3)$ & $ + 3n + 3$ &      &\ \ \ \ \ $\checkmark$ \\
\hline
$\SO(m_{1})\times\SO(m_{2})$ &    &   & \\
$m_{1} + m_{2} <$  &   &  & \\
$k_1 + k_2 + k_3 -1$ &   &  & \\
$n = (k_1 + k_2 + k_3)-$   & ${s = \binom{n}{2}}+$  &  & \\
$ (m_{1} + m_{2})$  & $ + 2n + 1$  &  & \ \ \ \ \ $\checkmark$\\
\hline
$\SO(d)$   &   &    &  \\
$d < k_1 + k_2 + k_3 -1$   &    &   &\\
$n = (k_1 + k_2 + k_3) - d$  & $s = \binom{n}{2} + n$  &     &\ \ \ \ \ $\checkmark$\\
\hline
$\SO(k_{3})$  &  &   &  \\
$n = k_{1} + k_{2}$  &  $s = \binom{n}{2} + n$ &    &\ \ \ \ \ $\checkmark$\\
\hline
\end{tabular}
\end{center}
Table 1: The number of isotropy summands for the homogeneous space $G/H = \SO(k_1 + k_2 + k_3)/H$.  The four last spaces contain  equivalent isotropy summands.
}
\end{example}
 
In the last four cases the complete description of $\Ad(H_{i})$-invariant inner products is much more difficult, because $\langle\fr{m}_{i}, \fr{m}_{j}\rangle$ are not necessarily zero for $i\neq j$.  
\begin{example}
\textnormal{We compute the complexified isotropy representation of the Stiefel manifold $V_{k}\bb{H}^{n}\cong \Sp(n)/\Sp(n-k)$, i.e. $\Ad^{\Sp(n)/\Sp(n-k)}\otimes\bb{C} : \Sp(n-k) \to \Aut(\fr{m}\otimes\bb{C})$.  It is
\begin{eqnarray*}
\left.\Ad^{\Sp(n)}\otimes\bb{C}\right|_{\Sp(n-k)} &=& S^{2}\nu_{n}\Big|_{\Sp(n-k)} = S^{2}(\nu_{n-k}\oplus1_{k}\oplus1_{k})\\
&=& S^{2}\nu_{n-k}\oplus S^{2}(1_{k}\oplus1_{k})\oplus(\nu_{n-k}\otimes(1_{k}\oplus1_{k})) \\
&=& S^{2}\nu_{n-k}\oplus\underbrace{1\oplus\cdots\oplus 1}_{\binom{2k+1}{2}}\oplus\underbrace{\nu_{n-k}\oplus\cdots\oplus \nu_{n-k}}_{2k}\\
&=& \Ad^{\Sp(n-k)}\otimes\bb{C} \oplus 1\oplus\cdots\oplus 1\oplus\nu_{n-k}\oplus\cdots\oplus \nu_{n-k},
\end{eqnarray*} 
so from Proposition \ref{isotrepr} we have that $\Ad^{\Sp(n)/\Sp(n-k)}\otimes\bb{C} = 1\oplus\cdots\oplus 1\oplus\nu_{n-k}\oplus\cdots\oplus \nu_{n-k}$.  Therefore the complexified tangent space $\fr{m}\otimes\bb{C}$ of $\Sp(n)/\Sp(n-k)$ can be written as a direct sum of $\binom{2k+1}{2}$ and $2k$ complex subspaces, of dimensions 1 and $2(n-k)$ respectively.  
}
\end{example}

\begin{example}
\textnormal{Consider the projective space $\bb{C}{\rm P}^{2n+1} \cong \Sp(n+1)/\Sp(n)\times\U(1)$.  Then according to Proposition \ref{isotrepr} the complexified isotropy representation of this space is determined by the equation
$$
\Ad^{\Sp(n+1)}\otimes\bb{C}\big|_{\Sp(n)\times\U(1)} = (\Ad^{\Sp(n)\times\U(1)}\otimes\bb{C})\oplus(\Ad^{\Sp(n+1)/\Sp(n)\times\U(1)}\otimes\bb{C}).
$$  
Observe that the dimension of the adjoint representation of $\Sp(n)\times\U(1)$ is $2n^{2} + n + 1$.  We now compute, 
\begin{eqnarray*}
\Ad^{\Sp(n+1)}\otimes\bb{C}\big|_{\Sp(n)\times\U(1)} &=& S^{2}\nu_{n+1}\big|_{\Sp(n)\times\U(1)} = S^{2}(\nu_{n}\oplus\mu_{1}\oplus\bar{\mu}_{1})\\
&=& S^{2}\nu_{n} \oplus S^{2}\mu_{1} \oplus S^{2}\bar{\mu}_{1}\\  &&\oplus(\nu_{n}\otimes\mu_{1})\oplus(\nu_{n}\otimes\bar{\mu}_{1})\oplus(\mu_{1}\oplus\bar{\mu}_{1})\\
&=& \big(S^{2}\nu_{n}\oplus(\mu_{1}\otimes\bar{\mu}_{1})\big)\oplus S^{2}\mu_{1}\oplus S^{2}\bar{\mu}_{1}\\
&&\oplus(\nu_{n}\otimes\mu_{1})\oplus(\nu_{n}\otimes\bar{\mu}_{1})\\
&=& \Ad^{\Sp(n)\times\U(1)}\otimes\bb{C} \oplus S^{2}\mu_{1}\oplus S^{2}\bar{\mu}_{1}\\
&&\oplus(\nu_{n}\otimes\mu_{1})\oplus(\nu_{n}\otimes\bar{\mu}_{1}),
\end{eqnarray*} 
where the fourth equality holds because the dimension of $S^{2}\nu_{n}$ $\oplus(\mu_{1}\otimes\bar{\mu}_{1})$ is equal to the dimension of the adjoint representation of $\Sp(n)\times\U(1)$.  Hence, the isotropy representation decomposes into a sum of four irreducible subrepresentations of dimensions $1, 1, 2n$ and $2n$ respectively, that is 
$$
\Ad^{\Sp(n+1)/\Sp(n)\times\U(1)}\otimes\bb{C} = S^{2}\mu_{1}\oplus S^{2}\bar{\mu}_{1}\oplus(\nu_{n}\otimes\mu_{1})\oplus(\nu_{n}\otimes\bar{\mu}_{1}).
$$  
Thus, the complexified tangent space $\fr{m}\otimes\bb{C}$ of $\Sp(n+1)/\Sp(n)\times\U(1)$ is written as a direct sum of four complex subspaces as $\fr{m}\otimes\bb{C} = \fr{p}_{1}\oplus\fr{p}_{2}\oplus\fr{p}_{3}\oplus\fr{p}_{4}$.  The real subspace $\fr{m}$ splits into two real subspaces of dimension $2$ and $4n$ respectively, i.e. $\fr{m} = \fr{m}_{1}\oplus\fr{m}_{2}$, where $\fr{m}\otimes\bb{C} = \fr{p}_{1}\oplus\fr{p}_{2}$ and $\fr{m}_{2}\otimes\bb{C} = \fr{p}_{3}\oplus\fr{p}_{4}$.} 
\end{example}

It is worth mentioning that W. Ziller in \cite{Zi} proved that the projective space $\bb{C}{\rm P}^{2n+1} \cong \Sp(n+1)/\Sp(n)\times\U(1)$ admits precisely two Einstein metrics. 

In general all Lie groups $G$ which act on the projective spaces $\bb{C}{\rm P}^{n}, \bb{H}{\rm P}^{n}$ and $\bb{C}_{\al}{\rm P}^{2}$ where classified by  Onishchik \cite{On}, according to the following table:
\begin{center}
\begin{tabular}{ l l l l }
\hline\hline
$G$ &                 $H$                              & $G/H$                      &  isotr. repr. \\
\hline\hline
$\SU(n+1)$         &     ${\rm S}(\U(1)\times \U(n))$           & $\mathbb{C}{\rm P}^{n}$    &     irreducible\\      
$\Sp(n+1)$         &        $\Sp(n)\times \Sp(1)$         & $\mathbb{H}{\rm P}^{n}$    &      irreducible\\    
$\F_{4}$           &     ${\rm Spin}(9)$                &       $\mathbb{C}_{a}{\rm P}^{2}$  &    irreducible\\
$\Sp(n+1)$         &     $\Sp(n)\times \U(1)$       &    $\mathbb{C}{\rm P}^{2n+1}$   & $\frak{m} = \frak{m}_{1}\oplus \frak{m}_{2}$\\
\hline
\end{tabular}
\end{center}
\centerline{Table 2: Transitive actions on projective spaces.}

\smallskip

Observe that in the first three cases the isotropy representations are irreducible, which  means that the only $G$-invariant metric on these spaces is the \textit{standard homogeneous Riemannian metric} (i.e. the metric induced by the negative of the Killing form $B$ of $\fr{g}$).  By J. Wolf (\cite{Wo}) this metric is Einstein.

\section{A special class of $G$-invariant metrics on $G/H$}
Let $G$ be a compact Lie group and $H$ a closed subgroup so that $G$ acts almost effectively on $G/H$.  Let $\fr{g},$ $\fr{h}$ be the Lie algebras of $G$ and $H$ and let $\fr{g} = \fr{h}\oplus\fr{m}$ be a reductive decomposition of $\fr{g}$ with respect to some $\Ad(G)$-invariant inner product on $\fr{g}$, i.e. $\Ad(h)\fr{m}\subset\fr{m}$ for all $h\in H$ where $\fr{m} \cong T_{{\it o}}(G/H)$, ${\it o} = e H$.  For $G$ semisimple, the negative of the Killing form $B$ of $\fr{g}$ is an $\Ad(G)$-invariant inner product on $\fr{g}$, therefore we can choose the above decomposition with respect to this form.  A Riemannian metric $g$ on a homogeneous space $G/H$ is called $G$-invariant if the diffeomorphism $\tau_{\al} : G/H \to G/H,$ $\tau_{\al}(g H) = \al g H$ is a isometry.  The following proposition gives a description of $G$-invariant metrics on homogeneous spaces.

\begin{prop}\label{oto}
Let $G/H$ be a homogeneous space.  Then there exists a one-to-one correspondence between:
\begin{enumerate}
\item $G$-invariant metrics $g$ on $G/H$
\item $\Ad^{G/H}$-invariant inner products $\left\langle \cdot, \cdot\right\rangle$ on $\fr{m}$, that is 
$$
\langle \Ad^{G/H}(h)X,\, \Ad^{G/H}(h)Y\rangle = \langle X,\, Y\rangle \quad \mbox{for all}\ X, Y\in \fr{m}, h\in H\ \mbox{and} 
$$
\item (if $H$ is compact and $\fr{m} = \fr{h}^{\perp}$ with respect to the negative of the Killing form $B$ of $G$) $\Ad^{G/H}$-equivariant, $B$-symmetric\footnote{Or $B(\cdot, \cdot)$-self-adjoint endomorphisms $\fr{m}$.} and positive definite operators $A : \fr{m}\to\fr{m}$ such that 
$$
\langle X, Y\rangle = B(A(X), Y).
$$
\end{enumerate}  
We call such an inner product $\Ad^{G}(H)$-invariant, or simply $\Ad(H)$-invariant 
\end{prop}
From the above proposition we can see that the set of all $\Ad(H)$-invariant inner products on $\fr{m}$ can be parametrized by $\Ad(H)$-equivariant, symmetric and positive definite operators $A : \fr{m}\to \fr{m}$.  Thus we have
\begin{equation*}
\mathcal{M}^{G} \longleftrightarrow \left\{A : \fr{m}\to \fr{m}
 \  \Big\vert \ 
 \begin{array}{l}
  \Ad(H)\mbox{-equivariant, symmetric}\\ 
  \mbox{and positive definite operator}
  \end{array}
   \right\}. 
 \end{equation*}
It is clear that if $\fr{m}$ decomposes into a direct sum of $\Ad(H)$-invariant irreducible and pairwise inequivalent modules $\fr{m}_{i}$ of dimension $d_{i}$ $(i = 1,\ldots,s)$, that is 
$
\fr{m} = \fr{m}_{1}\oplus\cdots\oplus\fr{m}_{s}, 
$
then all $\Ad(H)$-invariant inner products on $\fr{m}$ are given by
$$
\langle\cdot ,  \cdot\rangle = x_{1}(-B)|_{\fr{m}_{1}} + \cdots + x_{s}(-B)|_{\fr{m}_{s}}, \ x_{i}\in\bb{R}^{+}, \ i = 1,\ldots,s.
$$
In this case the matrix of the operator $A$ with respect to some $(-B)$-orthonormal adapted basis $\mathcal{B}$ of $\fr{m}$ is given by 
$$
[A]_{\mathcal{B}} = \begin{pmatrix}
x_1 {\rm Id}_{d_{1}} & &0  \\
 & \ddots &  \\
 0 & &  x_s {\rm Id}_{d_{s}}
\end{pmatrix}.
$$
In this case the $G$-invariant metrics are called \textit{diagonal}.  However, if the decomposition of $\fr{m}$ contains $r$ equivalent orthogonal modules $\fr{m}_{i}$, 
then the matrix of the operator $A$ with respect to some $(-B)$-orthonormal adapted basis $\mathcal{D}$ of $\fr{m}$ is given by 
$$
[A]_{\mathcal{D}} = \begin{pmatrix}
x_1 {\rm Id}_{d_{1}} & \al_{12} {\rm Id}_{d_{1}} & \cdots & \al_{1s} {\rm Id}_{d_{1}} \\
\al_{12} {\rm Id}_{d_{2}} & x_{2} {\rm Id}_{d_{2}} & \cdots & \al_{2s} {\rm Id}_{d_{2}}\\
\vdots & \vdots &\ddots & \vdots\\
\al_{1s}{\rm Id}_{d_{s}} & \al_{2s}{\rm Id}_{d_{s}} & \cdots & x_{s} {\rm Id}_{d_{s}}
\end{pmatrix}.
$$
The number of $\al_{ij}$ is $\frac{r(r-1)}{2}$. 
  For the unitary and symplectic case, we consider for simplicity the case where the decomposition of $\fr{m}$ contains two equivalent modules, say $\fr{m}_{1}\cong\fr{m}_{2}$, of dimension $d$.  Here there are two and four non diagonal elements respectively.  For example in the unitary case, the matrix of the operator $A$ with respect to some $(-B)$-orthonormal adapted basis $\mathcal{D}$ of $\fr{m}$ is given by linear combinations of the matrices  
$$
J_{1} = \begin{pmatrix}
0  & \al_{1}{\rm Id}_{2d}  \\
\al_{1}{\rm Id}_{2d}  &  0
\end{pmatrix}, \quad J_{2} = \begin{pmatrix}
0 & 0 &0 &  \al_{2}{\rm Id}_{d} \\
0 & 0 & -\al_{2}{\rm Id}_{d} & 0  \\
0 & -\al_{2}{\rm Id}_{d} & 0 & 0 \\
\al_{2}{\rm Id}_{d} & 0 & 0 & 0
\end{pmatrix}, 
$$
$\al_{1}, \al_{2}\in\bb{R}^{+}.$  The idea behind our approach is to try to eliminate some of the non diagonal elements in the above matrix, and restrict the study to the diagonal metrics.  For the same problem in the case of a Lie group, K. Y. Ha and J. B. Lee in \cite{HaLe}  classified the left-invariant Riemannian metrics for each simply connected three-dimensional Lie group \textit{up to automorphism}.  The main idea there was to identify all automorphisms of the Lie algebra of these groups, and then define an action of the automorphism group on the set of all left invariant inner products on the Lie algebras of these Lie groups\footnote{In general, the group of automorphisms of a Lie group $G$ defines an action on the set of all metrics on $G$ by $\Aut(G)\times \{\mbox{metrics on}\ G\} \to \{\mbox{metrics on}\ G\}$, $(\theta, g(\cdot, \cdot))\mapsto \theta g(\cdot, \cdot) := g_{\theta}(d\theta^{-1}\cdot, d\theta^{-1}\cdot)$.  Note that if the metric $g$ is left-invariant, then the metric $g_{\theta}$ is not necessarily left-invariant.}.  More precisely,
let $G$ be a Lie group and $\fr{g}$ the corresponding Lie algebra of $G$.  Let $\fr{M}$ be the set of all left invariant inner products of $\fr{g}$.  Then  $\Aut(\fr{g})$ acts on $\fr{M}$  by 
$$
\Aut(\fr{g})\times \fr{M} \to \fr{M}, \quad (\phi\, ,\, \left\langle \cdot, \cdot\right\rangle)\mapsto \left\langle \phi^{-1}\cdot\, ,\, \phi^{-1}\cdot\right\rangle.
$$
Under this action we can define an equivalence relation $\sim$ on $\fr{M}$ as follows:
$$
\langle \cdot, \cdot \rangle \sim \langle \cdot, \cdot\rangle'\ \Longleftrightarrow \ \mbox{there exists} \ \phi\in\Aut(\fr{g}) \ \mbox{such that} \ \langle \cdot, \cdot\rangle' = \langle \phi^{-1}\cdot\, ,\, \phi^{-1}\cdot\rangle. 
$$
Now, let $G/H$ be a homogeneous space ($H$ is the isotropy subgroup at the identity coset $e H$) with reductive decomposition $\fr{g} = \fr{h}\oplus\fr{m}$ with respect to some $\Ad(G)$-invariant inner product of $\fr{g}$.  Let $\Aut(G, H)$ be the set of all automorphisms of $G$ which preserve the group $H$.  It can be shown that if $\phi\in\Aut(G, H)$, then $\phi$ induces a $G$-equivariant diffeomorphism $\tilde{\phi} : G/H \to G/H$.  Then it is easy to see that this $G$-equivariant diffeomorphism defines an action on the set of all $G$-invariant metrics $\mathcal{M}^{G}$, transforming each $G$-invariant metric $g$ into a metric isometric to it.
  In general, every $G$-equivariant diffeomorphism of $G/H$ is a right translation by an element of $N_{G}(H)$, and for some $\al\in N_{G}(H)$ the map $\al \mapsto R_{\al}$, where $R_{\al} : G/H\to G/H$ is $G$-equivariant and sends each $gH$ to $g\al^{-1}H$.  This induces an isomorphism of $N_{G}(H)/H$ onto the group of $\Aut(G/H)$ (\cite{Br}).  Next, we describe when the set $\Aut(G/H)\cong N_{G}(H)/H$ defines an action on the set of all $G$-invariant metrics $\mathcal{M}^{G}$ of a homogeneous space $G/H$.   
  
First we recall the following fact.  Let $G_{1}$ and $G_{2}$ be Lie subgroups of a Lie group $G$.  If $G_{1}\subset G_{2}$, then $G_{1}$ is a subgroup of the Lie group $G_{2}$, and $\fr{g}_{1}\subset\fr{g}_{2}$.  Conversely, if $\fr{g}_{1}\subset\fr{g}_{2}$ and the group $G_{1}$ is connected, then $G_{1}\subset G_{2}$.  From this we have:
\begin{lemma}\textnormal{(cf. \cite{GOV})}\label{equal}
Let $G$ be a Lie group and $H$ be a closed, connected subgroup of $G$, with $\fr{g}$ and $\fr{h}$ the corresponding Lie algebras.  Then the group $N_{G}(H) = \{g\in G : gHg^{-1} = H\}$ is equal to the group $N_{G}(\fr{h}) = \{g\in G : \Ad(g)\fr{h}\subset\fr{h}\}$.
\end{lemma}
\begin{proof}
We need to show that (a) $N_{G}(H) \subset N_{G}(\fr{h})$ and (b) $N_{G}(\fr{h})\subset N_{G}(H)$.  For (a), let $g\in N_{G}(H)$.  Then $gHg^{-1} = H$ and by the above fact we have that $g\fr{h}g^{-1} = \fr{h}$, i.e. $\Ad(g)\fr{h} = \fr{h}$, hence $g\in N_{G}(\fr{h})$.  For (b), if $g\in N_{G}(\fr{h})$ then $\Ad(g)\fr{h}\subset\fr{h}$.  Since $\Ad(g)\fr{h}$ is the Lie algebra of $gHg^{-1}$ and $H$ is connected it follows that $gHg^{-1}\subset H$.  Obviously $H\subset gHg^{-1}$, hence we obtain that $g\in N_{G}(H)$. 
\end{proof}
The following proposition is central in our study.
\begin{prop}\label{protasi}
Let $n\in N_{G}(H)$ and  $\Ad(n): \fr{g}\to \fr{g}$.  Then the operator $\Ad(n)|_{\fr{m}} : \fr{m}\to \fr{g}$ takes values in $\fr{m}$, that is $\phi = \Ad(n)|_{\fr{m}}\in\Aut(\fr{m})$.  Also, $(\Ad(n)|_{\fr{m}})^{-1} = (\Ad(n)|_{\fr{m}})^{t}$.
\end{prop}
\begin{proof}
Let $n\in N_{G}(H)$ and $Y\in \fr{h}$.  Using Lemma \ref{equal}, for any subspace $\fr{h}$ of $\fr{g}$, the normalizer $N_{G}(\fr{h})$ is given by    
$
N_{G}(\fr{h}) = \{g\in G\ :\ \Ad(g)\fr{h}\subset \fr{h}\} = N_{G}(H).
$
Therefore, it follows that 
\begin{equation}\label{in}
\Ad(n)Y\in \fr{h}. 
\end{equation}
Let $X\in \fr{m} = \fr{h}^{\perp}$.  Then by using (\ref{in}) and the $\Ad(G)$-invariance of $B$ we obtain that
$$
B(\Ad(n)^{-1}X, Y) = B(\Ad(n)^{-1}X, \Ad(n)^{-1}\Ad(n)Y) = B(X, \Ad(n)Y) = 0,
$$ 
hence $\Ad(n)^{-1}X\in \fr{m}$.  Finally, for $n\in N_{G}(H)$ and using the $\Ad(G)$-invariance of $B$, we have that $B(\Ad(n)|_{\fr{m}}X,\, \Ad(n)|_{\fr{m}}Y) = B(X, Y)$. 
Since in general it is $B(\Ad(n)|_{\fr{m}}X,\, \Ad(n)|_{\fr{m}}Y)$ $= B(X, \, (\Ad(n)|_{\fr{m}})^{t}$ $\Ad(n)|_{\fr{m}}Y)$, it follows that $(\Ad(n)|_{\fr{m}})^{-1} = (\Ad(n)|_{\fr{m}})^{t}$.  
\end{proof}
Consider the set
$
\Phi = \{\phi = \Ad(n)|_{\fr{m}} : n\in N_{G}(H)\}.
$
Then by Proposition \ref{protasi} $\Phi$ is contained in $\Aut(\fr{m})$, hence we can define the isometric action\footnote{This action is essentially the action of $N_{G}(H)$ on $\mathcal{M}^{G}$, or equivalently the action of the group $N_{G}(H)/H$ on $\mathcal{M}^{G}$.} 
\begin{equation}\label{action}
\Phi \times \mathcal{M}^{G} \to \mathcal{M}^{G}, \quad (\phi\, ,\, A)\mapsto \phi\circ A\circ\phi^{-1}\equiv\tilde{A}.
\end{equation}
\begin{lemma}
The action of $\Phi$ on $\mathcal{M}^{G}$ is well defined.
\end{lemma}
\begin{proof}
We need to show that the operator $\tilde{A}$ is
\begin{itemize}
\item[(a)] $\Ad(H)$-equivariant, i.e.
$\Ad(H)\circ\tilde{A} = \tilde{A}\circ\Ad(H)$ or\\ $\Ad(H)\circ\tilde{A}\circ\Ad(H)^{-1} = \tilde{A}$ and
\item[(b)] $B$-symmetric and positive definite. 
\end{itemize}
For (a), let $n\in N_{G}(H)$ and we compute:
\begin{eqnarray*}
\Ad(H)\circ \tilde{A} \circ \Ad(H)^{-1} &=& \Ad(H)\circ \big(\Ad(n)\circ A \circ \Ad(n)^{-1}\big)\circ \Ad(H)^{-1}\\
&=& \Ad(Hn)\circ A \circ \Ad(Hn)^{-1} \\
&=& \Ad(nH)\circ A \circ \Ad(nH)^{-1}\\
&=& \Ad(n) \circ\big(\Ad(H)\circ A \circ \Ad(H)^{-1}\big) \circ \Ad(n)^{-1}\\
&=& \Ad(n)\circ A \circ \Ad(n)^{-1} = \tilde{A}.
\end{eqnarray*}
In the third equality we used the fact that $n\in N_{G}(H)$, 
and in the fifth equality the fact that the operator $A$ is $\Ad(H)$-equivariant.  For (b), let $X, Y\in \fr{m}$.  We will show that the operator $\tilde{A}$ is $B$-symmetric: 
\begin{eqnarray*}
B(\tilde{A}X, Y) &=& B(\Ad(n)\circ A \circ\Ad(n)^{-1}X, \, Y) \\
                 &=& B(\Ad(n)\circ A \circ\Ad(n)^{-1}X, \, \Ad(n)\Ad(n)^{-1}Y)\\
                 &=& B(A\circ \Ad(n)^{-1}X, \, \Ad(n)^{-1}Y) \\
                 &=& B(\Ad(n)^{-1}X,\,  A\circ\Ad(n)^{-1}Y)\\
                 &=& B(\Ad(n)\Ad(n)^{-1}X,\, \Ad(n)\circ A\circ\Ad(n)^{-1}Y)\\ 
                 &=& B(X, \tilde{A}Y).
\end{eqnarray*} 
In the third and fifth equality we used that the Killing form $B$ is $\Ad(G)$-invariant and in the fourth equality we used the fact that the operator $A$ is $B$-symmetric.  Finally, we show that $\tilde{A}$ is positive definite.  Let $X\in \frak{m}$ with $X\neq 0$.  Then by Proposition \ref{protasi} we have that $\Ad(n)X\in\fr{m}$ for all $n\in N_{G}(H)$, therefore it is
\begin{eqnarray*}
B(\tilde{A}X, X) &=& B(\Ad(n)\circ A \circ\Ad(n)^{-1}X, \, X) \\
                 &=& B(\Ad(n)\circ A \circ\Ad(n)^{-1}X, \, \Ad(n)\Ad(n)^{-1}X)\\
                 &=& B(A\circ \Ad(n)^{-1}X, \, \Ad(n)^{-1}X)\\
                 &=& B(A(\Ad(n)^{-1}X), \, \Ad(n)^{-1}X) > 0
\end{eqnarray*} 
where in the third equality we used that the Killing form $B$ is $\Ad(G)$-invariant.
\end{proof}
\begin{corol}
Let $n\in N_{G}(H)$.  Then the metrics corresponding to the operator $A$ are equivalent, up to the automorphism 
$
\Ad(n) : \fr{m} \to \fr{m},
$
to the metrics corresponding to the operator $\tilde{A}$. 
\end{corol} 

\begin{example}\textnormal{(\cite{Ke})}
\textnormal{We consider the Stiefel manifold $V_{2}\bb{R}^{4}\cong \SO(4)/\SO(2)$ and we will describe all $\SO(4)$-invariant metrics.  A reductive decomposition $\fr{m}$ of the  Lie algebra $\fr{so}(4)$, with respect to negative of Killing form $B(\cdot, \cdot)$ of $\SO(4)$, and for the embedding of
$
\fr{so}(2) \hookrightarrow 
\begin{pmatrix}
0 & 0\\
0 & \fr{so}(2)
\end{pmatrix}\in\fr{so}(4),
$
is the  set
\begin{eqnarray*}
\fr{m} &=& \left\{
\begin{pmatrix}
D_{2} & C\\
-C^{t} & O_{2}
\end{pmatrix} \, : \, D_{2} = {\rm diag}(0, 0),\, C\in M_{2}\bb{R}\right\} \\
&=& {\rm span}\{e_{ij} = E_{ij} - E_{ji} \, : \, 1\leq i < j \leq 4\},
\end{eqnarray*}
where $E_{ij}$ denotes the $4\times 4$ matrix with $1$ in $(ij)$-entry and $0$ elsewhere.  According to Example \ref{pisotr} the isotropy representation of $\SO(4)/\SO(2)$ is 
$\Ad^{\SO(4)/\SO(2)} = 1\oplus\lambda_{2}\oplus\lambda_{2}$,
thus $\fr{m}$ can be written as a direct sum of three $\Ad(\SO(2))$-invariant subspaces, of dimensions $1, 2, 2$ as follows
$$
\fr{m} = \fr{m}_{1}\oplus\fr{m}_{2}\oplus\fr{m}_{3}.
$$
Hence $\fr{m}_{1} = {\rm span}\{e_{12}\}, \, \fr{m}_{2} = {\rm span}\{e_{1j} \, : \,  j =3, 4\}$ and $\fr{m}_{3} = {\rm span}\{e_{2j} \, : \, j = 3, 4\}$.  Observe that $\fr{m}_{1} = \fr{so}(2) \hookrightarrow
\begin{pmatrix}
\fr{so}(2) & 0\\
0 & 0
\end{pmatrix}\in \fr{so}(4)$.  Since in the isotropy representation the  last two representations are equivalent and $\lambda_{2}\otimes\bb{C} \cong \lambda_{2}$, the space of intertwining maps is $1$-dimensional.
Therefore,  Proposition \ref{oto} implies that the matrix of the $\Ad(\SO(2))$-equivariant, symmetric and positive definite operator $A : \fr{m}\to \fr{m}$ with respect to some $(-B)$-orthonormal basis adapted to $\fr{m}$, is given by
$$
\left[A\right] = \begin{pmatrix}
x_{1}{\rm Id}_{1} & 0 & 0\\
0 & x_{2}{\rm Id}_{2} & \al{\rm Id}_{2}\\
0 & \al{\rm Id}_{2} & x_{3}{\rm Id}_{2}
\end{pmatrix} \quad \al\in\bb{R}^{+}.
$$
Here it is $N_{\SO(4)}(\SO(2)) = \SO(2)\cong S^{1}$.  The Lie algebra of $\SO(2)$ is generated by the element $e_{12} \hookrightarrow \fr{so}(4)$.  We consider the one parameter group $\exp(t e_{12})$ of $\SO(2)$, and we compute the matrix of the operator $\Ad(\exp(t e_{12})) : \fr{m} \to \fr{m}$ with respect to the basis $\{e_{ij}\}$ of $\fr{m}$.  We have the following:
\begin{eqnarray*}
\Ad(\exp(t e_{12}))e_{12} &=& e^{t e_{12}} e_{12} (e^{t e_{12}})^{-1} \\
&=& \begin{pmatrix}
\cos t & \sin t & 0 & 0\\
-\sin t & \cos t & 0 & 0\\
0 & 0 & 1 & 0\\
0 & 0 & 0 & 1
\end{pmatrix}\times
\begin{pmatrix}
0 & 1 & 0 & 0\\
-1 & 0 & 0 & 0\\
0 & 0 & 0 & 0\\
0 & 0 & 0 & 0
\end{pmatrix}
 \times\begin{pmatrix}
\cos t & -\sin t & 0 & 0\\
\sin t & \cos t & 0 & 0\\
0 & 0 & 1 & 0\\
0 & 0 & 0 & 1
\end{pmatrix} = e_{12}. 
\end{eqnarray*}  
Similarly, we obtain that
\begin{eqnarray*}
\Ad(\exp(t e_{12}))e_{13} &=& \cos t \cdot e_{13} -\sin t \cdot e_{23},\\
\Ad(\exp(t e_{12}))e_{14} &=&  \cos t \cdot e_{14} - \sin t \cdot e_{24},\\
\Ad(\exp(t e_{12}))e_{23} &=&  \sin t \cdot e_{13} + \cos t \cdot e_{23}, \\ 
\Ad(\exp(t e_{12}))e_{24} &=&  \sin t \cdot e_{14} + \cos t \cdot e_{24}.
\end{eqnarray*}
In the above calculations we used the fact that \textit{for any matrix group $G$ we have that $\Ad(g)X = g X g^{-1}$ for all $g\in G$, and $X\in \fr{g}$}.  Hence, the matrix of the operator $\Ad(\exp(t e_{12}))$ is
$$
\left[\Ad(\exp(t e_{12}))\right] = 
\begin{pmatrix}
1 & 0 & 0 & 0 & 0\\
0 & \cos t & 0 & \sin t & 0\\
0 & 0 & \cos t & 0 & \sin t\\
0 &  -\sin t & 0 & \cos t & 0\\
0 & 0 & -\sin t & 0 & \cos t 
\end{pmatrix}.
$$
If we set $\varphi = \Ad(\exp(t e_{12}))$, then the action (\ref{action}) at the matrix level is given by 
$$
\left(\left[\Ad(\exp(t e_{12}))\right], \left[A\right]\right)\mapsto \left[\Ad(\exp(t e_{12}))\right]\cdot\left[A\right]\cdot\left[\Ad(\exp(t e_{12}))\right]^{-1} \equiv [\tilde{A}].
$$
After some calculations we obtain that 
$$
[\tilde{A}] = 
\begin{pmatrix}
x_{1}{\rm Id}_{1} & 0 & 0\\
0 & k{\rm Id}_{2} & m{\rm Id}_{2}\\
0 & m{\rm Id}_{2} & c{\rm Id}_{2}
\end{pmatrix}, 
$$
where $k = x_{2} + (x_{3} - x_{2})\sin^{2}t + 2\al\sin t \cos t$, $m = (x_{3} - x_{2})\sin t\cos t +\al\cos 2t$ and $c = x_{3} + (x_{2} - x_{3})\sin^{2}t - 2\al\sin t\cos t$.  Obviously, we can find $t$ such that $m = (x_{3} - x_{2})\sin t\cos t +\al\cos 2t = 0$.  Therefore, without loss of generality, and since $A$ and $\tilde{A}$ are isometric, we can assume that the matrix of the operator $A$ is such that $\al = 0$ (i.e. diagonal).}
\end{example}
The important points in the previous example are that, we have exactly one non diagonal element in the matrix of the metric $\langle\cdot, \cdot\rangle = -B(A\cdot, \cdot)$ and that the normalizer of $\SO(2)$ in $\SO(4)$ is the circle $S^{1}\cong\SO(2)$.  These enable us to eliminate the non diagonal element, 
since it is possible to describe the complete action of the normalizer on the space of all $G$-invariant metrics.  This situation occurs in \cite{Ke}.  In the case where the group $N_{G}(H)$ (or $N_{G}(H)/H$) is isomorphic to some other Lie group (see for example \cite{Ni1}) then it is more complicated to
describe explicitly the action of $N_{G}(H)$ on $\mathcal{M}^{G}$, so we try to confine our study in a suitable subset of $\mathcal{M}^{G}$.  

From the action (\ref{action}) we obtain the following interesting consequenses.  Let     
$$
(\mathcal{M}^{G})^{\Phi} = \{A\in\mathcal{M}^{G} \, :\, \phi\circ A\circ \phi^{-1} = A \ \mbox{for all} \ \phi \in \Phi\}
$$
be the set of all fixed points\footnote{Let $G$ be a Lie group and $M$ a manifold.  Consider the action of $G$ on $M$:
$G \times M \to M, \ (g, m)\mapsto g\cdot m.$  The subset $M^{G} = \{m\in M \, : \, g\cdot m = m \ \mbox{for all} \ g\in G\}$ of $M$ is called the fixed point set of the action.} 
of the action $\Phi$ on $\mathcal{M}^{G}$ (which is subset of $\mathcal{M}^{G}$).  Any element of $(\mathcal{M}^{G})^{\Phi}$ parametrizes all $\Ad(N_{G}(H))$-invariant inner products of $\fr{m}$ and thus defines a subset of all inner product on $\fr{m}$.
Since $H\subset N_{G}(H)$, Proposition \ref{oto} can be restated as follows:
\begin{prop}\label{refre}
Let $G/H$ be a homogeneous space.  Then there exists a one-to-one correspondence between:
\begin{itemize}
\item[(1)] $G$-invariant metrics on $G/H$,
\item[(2)] $\Ad(H)$-invariant inner products $\langle\cdot, \cdot\rangle$ on $\fr{m}$,
\item[(3)] Fixed points $(\mathcal{M}^{G})^{\Phi_{H}} = \{A\in\mathcal{M}^{G} \,:\, \psi\circ A\circ \psi^{-1} = A, \ \mbox{for all} \ \psi \in \Phi_{H}\}$ of the action $\Phi_{H} = \{\phi = \Ad(h)|_{\fr{m}} : h\in H\}\subset \Phi$ on $\mathcal{M}^{G}$.
\end{itemize} 
\end{prop} 
Observe that $(\mathcal{M}^{G})^{\Phi}\subset(\mathcal{M}^{G})^{\Phi_{H}}$.  In the case where $N_{G}(H)\neq S^{1}$, we can work with some appropriate closed subset $K$ of the Lie group $G$, such that $H\subset K\subset N_{G}(H)$.  Then the fixed point set of the non trivial action of the set $\Phi_{K} = \{\varphi = \Ad(k)|_{\fr{m}} : k\in K\}\subset \Phi$ on $\mathcal{M}^{G}$ is
$$
(\mathcal{M}^{G})^{\Phi_{K}} = \{A\in\mathcal{M}^{G}\,:\, \varphi\circ A\circ \varphi^{-1} = A \ \mbox{for all} \ \varphi \in \Phi_{K}\},
$$
and this set determines a subset of all $\Ad(K)$-invariant inner products of $\fr{m}$.  We have the inclusions $(\mathcal{M}^{G})^{\Phi}\subset(\mathcal{M}^{G})^{\Phi_{K}}\subset(\mathcal{M}^{G})^{\Phi_{H}}$. 

\begin{figure}
\centering
\includegraphics[width=8cm]{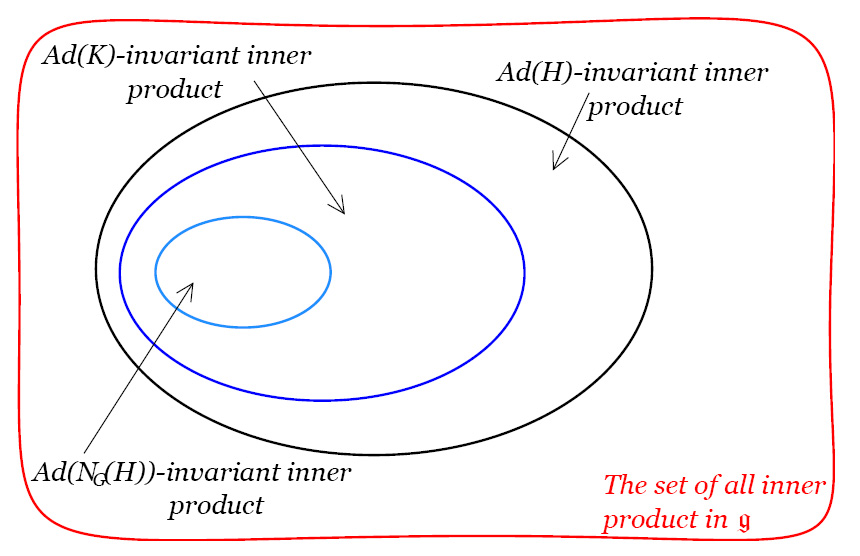}
\caption{Inclusions of certain invariant inner products in $\fr{g}$}
\end{figure}

By Proposition \ref{refre} the subset $(\mathcal{M}^{G})^{\Phi_{K}}$ is in one-to-one correspondence with a subset $\mathcal{M}^{G,K}$ of all $G$-invariant metrics, call it $\Ad(K)$-invariant, as shown in the Figure 2.

\begin{figure}
\centering
\includegraphics[width=9cm]{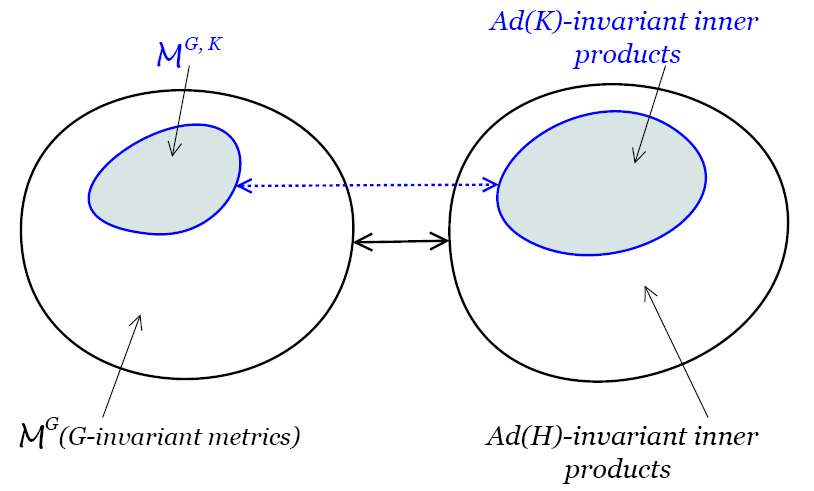}
\caption{Correspondence between $\Ad(K)$-invariant inner products on $\fr{g}$ and a subset of $G$-invariant metrics on $G/H$.}
\end{figure}

In the special case where $H = \{e\}$, then $N_{G}(H) = G$, thus the fixed points of the action (\ref{action}) are the $\Ad(G)$-invariant inner products on $\fr{g}$, which correspond to the bi-invariant metrics on the Lie group $G$.
 
We will now make  an appropriate choice of the subgroup $K$ in $G$.

\begin{prop}\label{subset}
Let $K$  be a subgroup of $G$ with $H\subset K \subset G$ and such that $K = L\times H$, for some subgroup $L$ of $G$.  Then $K$ is contained in $N_{G}(H)$.
\end{prop}
\begin{proof}  
We will show that if $k = (l , h)\in K = L\times H$ then $k H k^{-1} = H$ (that is $k\in N_{G}(H)$).
We identify $H$ with  $\{e\}\times H$, where $\{e\}$ is the identity element of $G$, 
so we will show that $k \big(\{e\}\times H \big) k^{-1} = \{e\}\times H$.  
It is
\begin{eqnarray*}
k \big(\{e\}\times H \big) k^{-1} &=& (l , h)\big(\{e\}\times H\big)(l, h)^{-1} = (l, h)\big(\{e\}\times H\big)(l^{-1}, h^{-1})\\
&=& (l, h)(e, X)(l^{-1}, h^{-1}), \ \mbox{for all} \ X\in H\\
&=& (l, h X)(l^{-1}, h^{-1}) = (l l^{-1}, h X h^{-1}) = (e, h X h^{-1}) \\
&=& \{e\}\times H,
\end{eqnarray*}   
hence $K = L\times H \subset N_{G}(\{e\}\times H) = N_{G}(H)$. 
\end{proof}
Also from \cite{Ro} we have the following result:
\begin{prop}
Let the group $G$ and the subgroup $H$ of $G$.  Then $H\triangleleft N_{G}(H)\leq G$, and whenever $H\triangleleft J\leq G$, then $J$ is a subgroup of $N_{G}(H)$ (here $A\triangleleft G$ means that $A$ is a normal subgroup of $G$).
\end{prop}

\section{$\Ad(K)$-invariant metrics and Riemannian submersions}
In the present section we will relate the $\Ad(K)$-invariant metrics on $G/H$ defined in the previous section, to Riemannian submersions.  For $H\subset K\subset G$ such that $K\subset N_G(H)$, we consider the fibration
$$
K/H\to G/H\to G/K.
$$
Let $\frak{a}$ and $\frak{p}$ be the orthogonal complements of $\fr{h}$ in $\fr{k}$ (i.e. $\fr{k} = \fr{h}\oplus\fr{a}$),  and of $\fr{k}$ in $\fr{g}$ (i.e. $\fr{g} = \fr{k}\oplus\fr{p}$), with respect to the negative of the Killing form of $\fr{g}$.  We assume that $\fr{a}$ is also $\Ad(K)$-invariant subspace of $\fr{k}$.  The spaces $\fr{a}$ and $\fr{p}$ are called \textit{vertical} and \textit{horizontal} subspaces of $\fr{g}$.  Then we have the decomposition $\fr{g}=\fr{h}\oplus\fr{m}=\fr{h}\oplus\fr{a}\oplus\fr{p}$.  

Any $\Ad(K)$-invariant inner product on $\fr{p}$ defines a $G$-invariant metric $\check{g}$ on $G/K$ and any
$\Ad(H)$-invariant inner product on $\fr{a}$ defines a $K$-invariant metric $\hat{g}$ on $K/H$.  The direct sum of these inner products on $\fr{a}\oplus\fr{p}$ defines a $G$-invariant metric
\begin{equation}\label{submersionmetric}
g=\hat{g}+\check{g}
\end{equation}
on $G/H$, called \textit{submersion metric}.  As this metric can be determined by an $\Ad(K)$-invariant inner product on $\fr{m}\cong T_{{\it o}}(G/H)=\fr{a}\oplus\fr{p}$,  it corresponds to an element of $(\mathcal{M}^G)^{\Phi_K}$, as defined in the previous section.  Hence we have the following:

\begin{prop}
Let $M=G/H$ be a homogeneous space and let $K$ be a closed subgroup of $G$ chosen as in the Proposition \ref{subset}.
Then the metric $(\ref{submersionmetric})$ is an element of $(\mathcal{M}^G)^{\Phi_K}$.
\end{prop}

\begin{example}
\textnormal{Let $M=G/H=\SO(k_1+k_2+k_3)/\SO(k_3)$ $(k_1, k_2, k_3\ge 2)$ be the Stiefel manifold $V_{k_1+k_2}\mathbb{R}^{k_1+k_2+k_3}$, and let $K=L\times H = (\SO(k_1)\times\SO(k_2))\times\SO(k_3)$.  Then by Proposition \ref{subset} it is $K \subset N_{G}(H)$.  We consider the fibration
$$
\xymatrix{
 \displaystyle{\frac{\SO(k_1)\times\SO(k_2)\times\SO(k_3)}{\SO(k_3)}} \ar[r] &
     \displaystyle{\frac{\SO(k_1 + k_2 + k_3)}{\SO(k_3)}} \ar[d] \\
&\displaystyle{\frac{\SO(k_1 + k_2 + k_3)}{\SO(k_1)\times\SO(k_2)\times\SO(k_3)}}}
$$ 
Then the base space $G/K$ is a generalized Wallach space and it is known by Example \ref{pisotr} (cf. (\ref{isotrwall})) 
that the isotropy representation is a direct sum of three non equivalent subrepresentations.  Therefore, the tangent space $\fr{p}\cong T_{{\it o}}(G/K)$ decomposes into three $\Ad(K)$-invariant and non equivalent modules   $\fr{p}_{12}\oplus\fr{p}_{13}\oplus\fr{p}_{23}$ of dimensions $k_{i}k_{j}, \ i\neq j$.  The tangent space $\fr{a}\cong T_{{\it o}}(K/H)$ of the fiber $K/H$ is the Lie algebra $\fr{so}(k_{1})\oplus\fr{so}(k_{2})$, and it is $\Ad(H)$-invariant.  Also $\fr{a}$ is $\Ad(K)$-invariant, thus the tangent space $\fr{m}\cong T_{{\it o}}(G/H)$ of the total space $G/H$ can be written as $\fr{m} = \fr{a}\oplus\fr{p}$, and thus decomposed into five $\Ad(K)$-invariant non equivalent modules: 
$$
\fr{m}=\fr{a}\oplus\fr{p}=\fr{so}(k_1)\oplus\fr{so}(k_2)\oplus\fr{p}_{12}\oplus\fr{p}_{13}\oplus\fr{p}_{23}.
$$ 
Therefore, any $\Ad(K)$-invariant metric is diagonal and determined by $\Ad(K)$-invariant inner products on $\fr{m}$ of the form:
\begin{equation*}
\begin{array}{lll}
\langle\cdot,\cdot\rangle &=&  x_1 \, (-B) |_{\fr{so}(k_1)}+ x_2 \, (-B) |_{ \fr{so}(k_2)}\\
&& + x_{12} \,  (-B) |_{ \fr{p}_{12}}+ x_{13} \,  (-B) |_{ \fr{p}_{13}} + x_{23} \,  (-B) |_{ \fr{p}_{23}}.  
 \end{array}
\end{equation*}
These inner products are the $\Ad(K)$-invariant inner products of Figure 2.
}
\end{example}
\noindent
New invariant Einstein metrics on the Stiefel manifold $V_{k_1+k_2}\mathbb{R}^{k_1+k_2+k_3} \cong$ $\SO (k_1+k_2+k_3)$ $/\SO(k_3)$, with respect to the above inner products, were studied in \cite{ArSaSt}.



\end{document}